\documentclass[12pt,twoside]{article}
\usepackage{smart}
\usepackage{emlines2}
\usepackage{ergomath}
\usepackage{amsopn}
\Aiv
\plainpage

\mathdef\cc#1$\cal{#1}$

\def\ar#1 #2,#3 #4,#5;{\put(#1,#2){\vector(#3,#4){#5}}}
\def\bx#1 #2,#3;{\put(#1,#2){\makebox[0pt]{$#3$}}}
\def\ln#1 #2,#3 #4,#5;{\put(#1,#2){\line(#3,#4){#5}}}
\def\mar#1 #2,#3 #4,#5;#6 #7,#8;{\multiput(#1,#2)%
(#6,#7){#8}{\vector(#3,#4){#5}}}
\def\mln#1 #2,#3 #4,#5;#6 #7,#8;{\multiput(#1,#2)%
(#6,#7){#8}{\line(#3,#4){#5}}}

\def\fnote#1{\footnote}

\mathdef\wt#1$\widetilde{#1}$
\mathdef\wh#1$\widehat{#1}$
\mathdef\ol#1$\overline{#1}$
\mathdef\ul#1$\underline{#1}$
\def\la{\langle}
\def\ra{\rangle}
\newtheorem{Th}{Theorem}[section]
\newtheorem{Prop}[Th]{Proposition}
\newtheorem{Cor}[Th]{Corollary}

\newtheorem{emrem}{Remark}[section]

\def\arr{\longrightarrow}
\def\harr{\hookrightarrow}

\mathdef\sar#1,#2,#3;$#1\stackrel{#2}{\arr}#3$
\newcommand{\rmnameii}[2]
  {\expandafter\newcommand \csname #1\endcsname {{\operatorname{#2}}}}
\newcommand{\bfname}[1]
  {\expandafter\newcommand \csname #1\endcsname {{\operatorname{\bf #1}}}}
\newcommand{\oplim}[1]
  {\expandafter\def\csname #1\endcsname {\operatornamewithlimits{#1}}}
\bfname{Groups}
\bfname{AbGroups}
\bfname{LocGroups}
\bfname{Germ}
\bfname{Marr}
\mathdef\grp#1$\Groups(\cc{#1})$
\mathdef\agrp#1$\AbGroups(\cc{#1})$
\mathdef\lgrp#1$\LocGroups(\cc{#1})$
\mathdef\germ#1$\Germ\,\cc{#1}$
\mathdef\marr#1$\Marr\,\cc{#1}$
\mathdef\U${\bf U}$
\mathdef\Ui${\bf U}'$
\mathdef\Uii${\bf U}''$

\begin{document}
\subjclass{18F99 (Primary), 58A05 (Secondary)}
\keywords{Differential calculus, category theory}
\title{Differential Calculus in Categories. I.}
\author{Vladimir Molotkov\thanks{This work was written with financial support of
MPIMiS (Leipzig) and
the Ministry of Sci. and Educ. of Bulgaria
grant F-610/98-99.}\\
Inst. for Nuclear Res. and Nucl. Energy\\
blvd.\ Tsarigradsko shosse 72, Sofia 1784, Bulgaria\\
{\sl e-mail: vmolot@inrne.bas.bg}}
\date{September 19, 2005}
\maketitle
\section{Introduction}
\label{intr}

We will deal here with the following general situation. Let \cc{C} be a
category equipped with both some ``topological'' structure consisting of some
subclass \cc{O} of mono's of \cc{C} (open mono's) and with some ``algebraic''
structure consisting of a faithful functor
$G\colon\cc{V}\arr\grp C$, where \grp C
is the category of group objects in \cc{C}. For example: \cc{C} is the
category {\bf Top} of topological spaces, \cc{O} is the set of all mono's
of \cc{C}, isomorphic to inclusions of open subspaces, \cc{V} is the category
of all finite-dimensional (or Banach, or locally convex) spaces over
the field of real numbers, whereas $G$ is the forgetful functor.

The problem to be studied here is the following one: what means, in the
situation above, the approximation of arrows (in \cc{C}) between
``open subobjects of objects of \cc{V}'' by arrows which ``locally'' are
morphisms of \cc{V}? A part of this problem is: how to interpret such thing
as the first difference $\Delta f=f(u+v)-f(u)$, arising in classical
calculus, in a more general situation, when objects of the category \cc{C}
are not determined by their ``points''? Such ``non-classical'' situations
really exist in nature: for example, the theory of smooth supermanifolds of Berezin
and Leites \cite{BeLe} (which originates from
its algebraic geometric version~\cite{Le}),
was reformulated by the author (see~\cite{Mol}) as a
theory of ``smooth objects'' in some functor category, which permitted one
to generalize it to infinite dimensions. In fact, just the latter observation
has initiated the present author's work\footnote{In fact,
the essential part of the present work
was written in 1990. The author delayed its publication
in the hope to find a time to develop the theory a bit further, including
Taylor series expansions in it. This is the first serious "branching point":
it turned out to be not so easy to find "the" true categorical generalization
of both Taylor series and of manifolds of class $C^k$ for $k>1$. Of course,
the standard iterative definition of manifolds of class $C^k$ works, but
one needs to find further {\it categorical} restrictions on the calculus
in order that the Taylor series expansion were valid for morphisms of
class $C^k$ for $k>1$.

Besides, the author hoped to find a suitable extension
of ``classic'' calculus on locally convex vector spaces to a
much wider class of ``local models'', including, in particular,
arbitrary products of finite-dimensional Lie groups, both abelian
and non-abelian. This includes both the search of an adequate
generalization of the very notion of local convexity to the case
of topological groups, generated by their 1-parameter subgroups
and the class of ``small'' continuous maps between such
locally convex topological groups.
Now this generalization is found, at last, giving a quite new
example for the theory developed. Detailed construction
of calculus on locally convex groups will be given elsewhere.
%
}.

It turns out that the ``topology'' \cc{O} on the category \cc{C} generates
the category of germs of neighborhoods of mono's \germ C
objects of which are mono's of \cc{C} taken together with their
``infinitesimal neighborhood''
(see sect. \ref{germ} below for exact definition); besides, it generates a
functor attaching to a group $G$ of \cc{C} some group $\wt{\cc{C}}$
of the category \germ C which is {\sl local} in the sense
that it has the only point (namely, the unit ``element'' of the group).
At last, any arrow $f\colon U\arr U'$ between open subobjects of groups
$G$, resp. $G'$ generates in \germ C the {\sl difference arrow}
$\Delta f\colon U\times \wt{G}\arr \wt{G'}$
(sect. \ref{difmor}); it is just the latter arrow that is afterwards to be
approximated by choosing two classes \cc{L} (``linear'') and $o$ (``small'')
of  ``families'' of arrows having trivial intersection. The data \cc{C},
\cc{V}, \cc{L} and $o$ are called {\sl differential calculus} if they satisfy
some natural conditions guaranteeing both the functorial construction of the
category of smooth objects and of the ``tangent functor'' (sect. \ref{calc}).

Note that we are dealing here with calculus which is, generally speaking,
non-abelian; none of constructions of the present paper require that groups
under consideration are abelian. But an essential problem is not considered:
what abstract conditions one needs to impose on differential calculus in
order that something like Taylor expansion were valid for smooth arrows.
In other words, what is the categorical counterpart of mean-value theorem or,
rather, of its weakest corollaries, permitting one to deduce both Taylor
formula and relation between derivative and partial derivatives. The author
hopes to return to this question elsewhere.

\section{Locuses}
\label{loc}
A category \cc{C} together with some class \cc{O} of pullbackable mono's of
\cc{C} will be called a {\bf locus}, if \cc{O} contains all iso's and is
closed both under pullbacks and under the composition of arrows. Elements of
\cc{O} will be called {\bf open arrows} or, simply, {\bf opens}.

For example, if $\tau$ is some pretopology on \cc{C}such that every element of
every covering of $\tau$ is mono, then one can take
$\cc{O}=\cup{\rm Cov}\tau$.

Given some arrow $f\colon C\arr X$ of \cc{C} an open arrow
$U\stackrel{u}{\arr}X$ will be called a {\bf neighborhood} or, briefly,
a {\bf nood} of $f$ if there exists $f'\colon C\arr U$ such that $f=uf'$;
if $g\colon X\arr Y$ is another arrow then the arrow $gu$ will be called the
{\bf restriction} of $g$ on $u$ (or on $U$) and will be often denoted as
$g|_U$.

If $(\cc{C},\cc{O})$ and $(\cc{C}',\cc{O}')$ are two locuses, then a functor
$F\colon \cc{C}\arr\cc{C}'$ will be called a {\bf morphism of locuses} if
$F(\cc{O})\subset\cc{O}'$ and $F$ respects pullbacks of opens along
arbitrary change of base.

\section{Germs of noods of mono's}
\label{germ}
From now on \cc{C} will be some fixed locus. Let \marr C be
the full subcategory of the category of commutative squares of \cc{C}, whose
objects are just all mono's of \cc{C}. Let $S$ be the class of all morphisms
of \marr C of the form
\begin{equation}
\unitlength=0.50mm
\special{em:linewidth 0.4pt}
\linethickness{0.4pt}
\begin{picture}(65.00,45.00)
\put(3.00,3.00){\makebox(0,0)[cc]{$U$}}
\put(63.00,3.00){\makebox(0,0)[cc]{$X$}}
\put(63.00,43.00){\makebox(0,0)[cc]{$C$}}
\put(3.00,43.00){\makebox(0,0)[cc]{$C$}}
\put(33.00,45.00){\makebox(0,0)[cb]{${\rm Id}$}}
\put(33.00,5.00){\makebox(0,0)[cb]{$u$}}
\put(7.00,43.00){\vector(1,0){52.00}}
\put(7.00,3.00){\vector(1,0){52.00}}
\put(63.00,36.00){\vector(0,-1){29.00}}
\put(3.00,36.00){\vector(0,-1){29.00}}
\emline{3.00}{36.00}{1}{5.00}{38.00}{2}
\emline{3.00}{36.00}{3}{1.00}{38.00}{4}
\emline{63.00}{36.00}{5}{65.00}{38.00}{6}
\emline{63.00}{36.00}{7}{61.00}{38.00}{8}
\end{picture}
\end{equation}
where $u$ is open.

\begin{Prop}
The class $S$ admits the calculus of right fractions {\rm(see, e.g.,
Ref.\cite{GaZi} for definitions).}
\end{Prop}

The category of fractions $\marr C[S^{-1}]$ will be called
the {\bf category of germs of nood's of mono's} of the locus \cc{C} and
will be denoted as \germ C.

It can be described as follows. Objects of \germ C are all
mono's of \cc{C}; a morphism $c\colon C\arr X$ into
$c_1\colon C_1\arr X_1$ is the equivalence class of the commutative diagram
\begin{equation}
\unitlength=0.50mm
\special{em:linewidth 0.4pt}
\linethickness{0.4pt}
\begin{picture}(86.00,68.00)
\put(4.00,2.00){\makebox(0,0)[cc]{$X$}}
\put(84.00,2.00){\makebox(0,0)[cc]{$X_1$}}
\put(84.00,62.00){\makebox(0,0)[cc]{$C_1$}}
\put(4.00,62.00){\makebox(0,0)[cc]{$C$}}
\put(24.00,32.00){\makebox(0,0)[cc]{$U$}}
\put(28.00,30.00){\vector(2,-1){52.00}}
\put(84.00,58.00){\vector(0,-1){52.00}}
\put(4.00,58.00){\vector(0,-1){52.00}}
\put(8.00,62.00){\vector(1,0){72.00}}
\put(6.00,59.00){\vector(2,-3){16.00}}
\put(22.00,29.00){\vector(-2,-3){16.00}}
\put(44.00,64.00){\makebox(0,0)[cb]{$f$}}
\put(54.00,19.00){\makebox(0,0)[lb]{$g$}}
\put(17.00,47.00){\makebox(0,0)[lb]{$d$}}
\put(17.00,17.00){\makebox(0,0)[lt]{$u$}}
\put(2.00,32.00){\makebox(0,0)[rc]{$c$}}
\put(86.00,32.00){\makebox(0,0)[lc]{$c_1$}}
\end{picture}
\end{equation}
where $u$ is a nood of $c$, with respect to the following equivalence
relation: given another such diagram $(2')$ with data
$(d',u',f',g')$, the diagrams (2) and $(2')$ are equivalent if there
exists a nood $V\stackrel{v}{\arr}X$ of $c$,
such that $v$ is a restriction of both $u$ and $u'$
(i.e. factors through their intersection)
and $g|_V=g'|_V$.

\begin{Prop}
\label{uniq}
If the diagrams $(2)$ and $(2')$ belong to the same equivalence class {\rm
(with respect to the equivalence relation just introduced)} then $f=f'$.
\end{Prop}

The description of the composition of arrows in \germ C is
left to the reader.

Define now the functors
\begin{equation}
\label{UandF}
U\colon\cc{C}\arr\germ C\ ,\ \ F\colon\germ C\arr\cc{C}
\end{equation}
as follows: $U(C\stackrel{f}{\arr}C_1)$ is an equivalence class of the diagram
\begin{equation}
\unitlength=0.50mm
\special{em:linewidth 0.4pt}
\linethickness{0.4pt}
\begin{picture}(65.00,49.00)
\put(7.00,43.00){\vector(1,0){52.00}}
\put(7.00,3.00){\vector(1,0){52.00}}
\put(63.00,36.00){\vector(0,-1){29.00}}
\put(3.00,36.00){\vector(0,-1){29.00}}
\emline{3.00}{36.00}{1}{5.00}{38.00}{2}
\emline{3.00}{36.00}{3}{1.00}{38.00}{4}
\emline{63.00}{36.00}{5}{65.00}{38.00}{6}
\emline{63.00}{36.00}{7}{61.00}{38.00}{8}
\put(63.00,43.00){\makebox(0,0)[cc]{$C'$}}
\put(3.00,43.00){\makebox(0,0)[cc]{$C$}}
\put(3.00,3.00){\makebox(0,0)[cc]{$C$}}
\put(63.00,3.00){\makebox(0,0)[cc]{$C'$}}
\put(65.00,22.00){\makebox(0,0)[lc]{${\rm Id}$}}
\put(1.00,22.00){\makebox(0,0)[rc]{${\rm Id}$}}
\put(33.00,45.00){\makebox(0,0)[cb]{$f$}}
\put(33.00,5.00){\makebox(0,0)[cb]{$f$}}
\end{picture}\ ,
\end{equation}
whereas for a
morphism $\gamma\colon c\arr c_1$ in \germ C represented by
the diagram (2), define $F(\gamma)$ as $f$ (this definition is correct due to
Prop.\ref{uniq}).

\begin{Prop}
\label{funs}
The functor $U$ is fully faithful; the functor $F$ is left adjoint to $U$.
\end{Prop}

Noting that $FU={\rm Id}_{\cc{C}}$ one can {\it identify} \cc{C} with a full
subcategory of \germ C; due to Prop.\ref{funs} the category
\cc{C} is a reflective subcategory of \germ C.

Let $\cc{O}'$ be the set of all morphisms of \germ C having a
representation (2) such that both $f$ and $g$ are open and, besides, the square
of the diagram (2) is a pullback.

\begin{Prop}
The pair $(\germ C,\cc{O}')$ is a locus; moreover, both $F$ and
$U$ are morphisms of locuses.
\end{Prop}

\begin{Prop}
\label{prod}
If \cc{C} has finite products then \germ C has them as well,
whereas the functors $F$ and $U$ respect finite products. The product of
objects $f\colon C\arr X$ and $f'\colon C'\arr X'$ is $f\times f'$ with
projection arrows represented by the diagram
\begin{equation}
\unitlength=0.5mm
\special{em:linewidth 0.4pt}
\linethickness{0.4pt}
\begin{picture}(102.00,48.00)
\put(0.00,0.00){\makebox(0,0)[cc]{$X$}}
\put(50.00,0.00){\makebox(0,0)[cc]{$X\times X'$}}
\put(100.00,0.00){\makebox(0,0)[cc]{$X'$}}
\put(100.00,40.00){\makebox(0,0)[cc]{$C'$}}
\put(50.00,40.00){\makebox(0,0)[cc]{$C\times C'$}}
\put(0.00,40.00){\makebox(0,0)[cc]{$C$}}
\put(0.00,36.00){\vector(0,-1){32.00}}
\put(50.00,36.00){\vector(0,-1){32.00}}
\put(100.00,36.00){\vector(0,-1){32.00}}
\put(62.00,40.00){\vector(1,0){32.00}}
\put(64.00,00){\vector(1,0){32.00}}
\put(36.00,00){\vector(-1,0){32.00}}
\put(36.00,40.00){\vector(-1,0){32.00}}
\put(-2.00,20.00){\makebox(0,0)[rc]{$f$}}
\put(52.00,20.00){\makebox(0,0)[lc]{$f\times f'$}}
\put(102.00,20.00){\makebox(0,0)[lc]{$f'$}}
\put(79.00,44.00){\makebox(0,0)[cb]{${\rm pr}_{C'}$}}
\put(23.00,44.00){\makebox(0,0)[cb]{${\rm pr}_C$}}
\put(23.00,4.00){\makebox(0,0)[cb]{${\rm pr}_X$}}
\put(79.00,4.00){\makebox(0,0)[cb]{${\rm pr}_{X'}$}}
\end{picture}
\end{equation}

If, in addition, both $f$ and $f'$ are open in \cc{C} then $f\times f'$
is open as well.
\end{Prop}


\section{Groups and group bundles in categories}
\label{gr}
Definitions of groups and another kinds of algebraic objects
in categories with products (see e.g., \cite{Sch}) will be slightly
reformulated here in a form dealing with
algebraic relations rather than with commutative diagrams. This will permit
to check in many cases the commutativity of diagrams arising by purely
algebraic calculations just in the same way as it is processed in the
case of ordinary groups, group actions, etc. in the category of sets.

Note as well that if one uses for projection arrows like
${\rm pr}_{X_\alpha}\colon X_1\times\dots\times X_n\arr X_\alpha$ more
``human'' notations as $x_\alpha$ (so that
$(x_1,\dots,x_n)={\rm Id}_{X_1\times\dots\times X_n}$) then algebra in a
category \cc{C} with finite products will not be distinguished from that
in the category of sets at least notationally; the difference will arise
only when one {\sl interprets} expressions like $(x_1,\dots,x_n)$: this is a
particular {\sl arrow} in the general case, and a general {\sl element} of the
product in the case of sets.

Let \cc{D} be a category with finite products. Given an object $G$ of
\cc{D} together with a triple of arrows
\begin{equation}
\label{grarr}
G\times G\stackrel{.}{\arr}G,\ \ G\stackrel{i}{\arr}G,\ \
p\stackrel{e}{\arr}G
\end{equation}
(here $p$ is a final object of \cc{D}), one can define for any object
$Z$ and any pair of arrows $Z\stackrel{f}{\arr}G$
and~$Z\stackrel{f'}{\arr}G$ the arrows $f\cdot f'$, $f^{-1}$ and~$1$ as follows:
\begin{equation}
\label{homgr}
f\cdot f':=(Z\stackrel{(f,f')}{\arr}G\times G\stackrel{.}{\arr}G);\
f^{-1}:=(Z\stackrel{f}{\arr}G\stackrel{i}{\arr}G);\
1:=(Z\arr p\stackrel{e}{\arr}G).
\end{equation}

One has, in particular:
\begin{equation}
\label{multeq}
\cdot=g_1\cdot g_2,\ \ i=g^{-1},
\end{equation}
where $(g_1,g_2)={\rm Id}_{G\times G}$, $g={\rm Id}_G$.

One trivially checks that for any $f,f'\colon Z\arr G$ and
$u\colon U\arr Z$ the following identities hold:
\begin{userlabel}
\begin{equation}
\label{arr-mult}
(f\cdot f')u=fu\cdot f'u,\ f^{-1}u=(fu)^{-1},\ 1u=1.
\end{equation}

Here and in what follows it is assumed implicitly that categorical
composition of arrows (denoted simply by concatenation of arrows) bounds
closer than multiplication $\cdot$ so that, for example, $fu\cdot f'u$
means always the same as the expression $(fu)\cdot(f'u)$ and not, say,
$f(u\cdot f')u$, even if the latter expression makes sense.

Note as well that we will always identify, for any category \cc{C},
the unary product functor $P\colon\cc{C}\arr\cc{C}$ with identity functor
${\rm Id}_{\cc{C}}$, so that for any arrow $x$ of any category $\cc{C}$
the identities $x=(x)=((x))=\dots$ hold; in particular, in the situation
$\sar Z,x,X;\sar{},f,Y;$ one always can write $f(x)$ in place of the
compposition arrow $fx$, to make algebraic relations in arbitrary category
\cc{C} (expressing commutativity of some diagram in \cc{C})
to look {\it exactly\/} as the corresponding relations between
elements of algebra of some kind in the category of sets. For example,
the equalities~\ref{arr-mult} can be rewritten as well as
follows:

\begin{equation}
\label{}
(f\cdot f')(u)=f(u)\cdot f'(u)\ f^{-1}(u)=(f(u))^{-1},\ 1(u)=1.
\end{equation}
\end{userlabel}

An object $G$ of \cc{D} together with a triple of arrows (\ref{grarr}) is
called a {\bf group} in \cc{D} if the following identities are valid:
\begin{equation}
\label{grdef}
(g_1\cdot g_2)\cdot g_3=g_1\cdot(g_2\cdot g_3),\ g\cdot g^{-1}=g^{-1}\cdot g=1,
\ 1\cdot g=g\cdot1=g,
\end{equation}
where $(g_1,g_2,g_3)={\rm Id}_{G\times G\times G}$ and $g={\rm Id}_G$.
It then follows that the identities (\ref{grdef}) are valid for any triple
of arrows $g_1,g_2,g_3\colon Z\arr G$, resp. for any arrow $g\colon Z\arr G$,
and not only for "generic" projection arrows $g_1,g_2,g_3$ and~$g$ of
(\ref{grdef}). In other words, for any object $Z$ of \cc{D} the hom-set
$[Z,G]$ is a group with respect to operations (\ref{homgr}).
Similar is true for another kinds of algebraic structures: a relation between
some algebraic operations is universally valid if it is valid for certain
generic projection arrows.

The group $G$ is called
{\bf abelian} if the identity
\begin{equation}
g_1\cdot g_2=g_2\cdot g_1
\end{equation}
is valid for $(g_1,g_2)={\rm Id}_{G\times G}$. For abelian groups we will
write further $f_1+f_2$, resp. $-f$, resp. $0$ instead of $f_1\cdot f_2$,
resp. $f^{-1}$, resp. $1$.

If $G$ is a group in \cc{C} we will denote \sar G\times\dots\times G,\cdot,G;
the arrow of $n$-fold multiplication for any natural number $n\neq 0$; for
this arrow the identity
\begin{equation}
\cdot=g_1\cdot\dots\cdot g_n\ \ \
((g_1,\dots,g_n)={\rm Id}_{G\times\dots\times G})
\end{equation}
holds, generalizing the identity (\ref{multeq}).

Now, given groups $G$ and $G'$ in \cc{D}, an arrow $f\colon G\arr G'$ of \cc{D}
is called a {\bf morphism of groups} if the identities
\begin{equation}
\label{grmordef}
f(g_1\cdot g_2)=fg_1\cdot fg_2,\ f(g^{-1})=(fg)^{-1},\
f1=1
\end{equation}
hold, where $(g_1,g_2)={\rm Id}_{G\times G}$ and $g={\rm Id}_G$. It then
follows that equalities (\ref{grmordef}) are valid for any triple of arrows
$g_1,g_2,g\colon Z\arr G$. In other words, for any object $Z$ of \cc{D}
the map $[Z,f]\colon[Z,G]\arr[Z,G']$ of hom-sets is morphism of groups.
Similar is true for another kinds of algebraic structures.

More generally, given an object $U$ of \cc{D}, an arrow
$f\colon U\times G\arr G'$ is called an $U$-{\bf family of morphisms
of groups} if the identities
\begin{equation}
f(u,g_1\cdot g_2)=f(u,g_1)\cdot f(u,g_2),\
f(u,g^{-1})=(f(u,g))^{-1},\ f(u,1)=1
\end{equation}
hold, where $(u,g_1,g_2)={\rm Id}_{U\times G\times G}$, resp.
$(u,g)={\rm Id}_{U\times G}$, resp. $u={\rm Id}_U$.

Groups in \cc{D} together with morphisms of groups form a category
which will be denoted as \grp D.

More generally, one can define the category ${\rm GBun}_{\times}\cc{D}$
of {\bf trivial group bundles} of the category \cc{D} as follows. Objects of
${\rm GBun}_{\times}\cc{D}$ are pairs $\la U,G\ra$ consisting of
an object $U$ of \cc{D} (the {\bf base} of the bundle) and a group $G$
of \cc{D} (the {\bf fiber} of the bundle), whereas a morphism
from $\la U,G\ra$ to $\la U',G'\ra$ is a pair of arrows
$\la\sar U,u,U';,\sar U\times G,f,U'\times G';\ra$ such that
the diagram
\begin{equation}
\unitlength=0.50mm
\special{em:linewidth 0.4pt}
\linethickness{0.4pt}
\begin{picture}(72.00,44.00)
\put(70.00,38.00){\vector(0,-1){32.00}}
\put(10.00,38.00){\vector(0,-1){32.00}}
\put(10.00,2.00){\makebox(0,0)[cc]{$U$}}
\put(70.00,2.00){\makebox(0,0)[cc]{$U'$}}
\put(70.00,42.00){\makebox(0,0)[cc]{$U'\times G'$}}
\put(10.00,42.00){\makebox(0,0)[cc]{$U\times G$}}
\put(72.00,22.00){\makebox(0,0)[lc]{${\rm pr}_{U'}$}}
\put(8.00,22.00){\makebox(0,0)[rc]{${\rm pr}_U$}}
\put(40.00,44.00){\makebox(0,0)[cb]{$f$}}
\put(40.00,4.00){\makebox(0,0)[cb]{$u$}}
\put(24.00,42.00){\vector(1,0){32.00}}
\put(14.00,2.00){\vector(1,0){52.00}}
\end{picture}
\end{equation}
is commutative and, besides, the composition arrow
\begin{equation}
\sar U\times G,f,U'\times G';\sar{},{\rm pr}_{G'},G';
\end{equation}
(the {\bf principal part} of the morphism $\la u,f\ra$)
is an $U$-family of morphisms of groups. The category \grp D
can be identified, clearly, with the full subcategory of
${\rm GBun}_{\times}\cc{D}$ consisting of group bundles with base being
the final object $p$. Note that any arrow $u\colon U\arr U'$ together with
any $U$-family $f\colon U\times G\arr G'$ of group morphisms determine the
morphism
\[
\la u,(u{\rm pr}_U,f)\ra\colon\la U,G\ra\arr
\la U',G'\ra
\]
of group bundles with principal part $f$.

Similarly, one can define actions of groups in a category on objects
of a category (see Sect.\ref{act} below) as well as rings and modules in
a category, etc.. For every kind {\bf T} of algebraic structure thus defined,
one can define the notion of a family of morphisms of {\bf T}-algebras and the
corresponding category ${\bf T}{\rm Bun}_{\times}\cc{D}$ of
{\bf trivial bundles of T-algebras}.

For the initial theory {\bf T} (with empty set of operations) we will write
${\rm Bun}_{\times}\cc{D}$
instead of ${\bf T}{\rm Bun}_{\times}\cc{D}$ and will call objects of the
latter category {\bf trivial bundles}. Denote as well
$.{\rm Bun}_{\times}\cc{D}$ the category of {\bf pointed trivial bundles}
arising from the theory {\bf T} with the only operation of arity zero (and
with no other operations).

\noindent{\bf Remark.} If \cc{C} is a locus with finite products such that
its set \cc{O} of opens arouse from some pretopology $\tau$
(see Sect.\ref{loc}) then for any algebraic theory {\bf T} one can define the
category ${\bf T}{\rm Bun}\cc{C}$ of {\bf locally trivial bundles of
T-algebras} of the category (or, rather, presite) \cc{C}. This is left
to the reader.

\section{Local groups in the category \germ C}
\label{locgr}
In this section \cc{C} will be some locus with finite products; $p$ will be
some fixed final object in \cc{C}, whereas products are supposed to be chosen
in such a way that for any object $Z$ of \cc{C} one has
$p\times Z=Z\times p=Z$ and the corresponding projection arrows are equal to
${\rm Id}_Z$ (this choice is always possible).

Arrows of the kind $p\stackrel{x}{\arr}X$ will be called {\bf points} (of $X$);
they are mono, hence, every point is an object of \germ C;
a point, considered as an object of  \germ C will be called
a {\bf point-like object}.

Note that ${\bf p}=(p\stackrel{{\rm Id}}{\arr}p)$ is a final object in
\germ C and it again is canonical, i.e.
${\bf p}\times {\bf Z}={\bf Z}\times {\bf p}={\bf Z}$ for every object
${\bf Z}=(Z\arr X)$ of \germ C.

A group $G$ in \germ C which is, simultaneously, a point-like
object will be called a {\bf local group} (N.Bourbaki would prefer the term
``groupuscule'').

Now we will show that there are many {\sl local} groups in the category
\germ C  (provided there are many groups in the category \cc{C}
itself).

Let $G$ be a group in \cc{C}, with multiplication
$G\times G\stackrel{.}{\arr}G$, unity $G\stackrel{i}{\arr}G$ and inverse
$p\stackrel{e}{\arr}G$. It is as well a group in \germ C
due to Prop.\ref{prod} and the remark after Prop.\ref{funs}. Let
$\wt G$ denotes the arrow $p\stackrel{e}{\arr}G$, considered as an
object of \germ C. We have the canonical mono
\begin{equation}
\label{locincl}
\iota_G\colon\wt G\arr G;
\end{equation}
in \germ C, represented by the diagram
\begin{equation}
\unitlength=0.50mm
\special{em:linewidth 0.4pt}
\linethickness{0.4pt}
\begin{picture}(64.00,42.00)
\put(6.00,40.00){\vector(1,0){52.00}}
\put(6.00,0.00){\vector(1,0){52.00}}
\put(62.00,33.00){\vector(0,-1){29.00}}
\put(2.00,33.00){\vector(0,-1){29.00}}
\emline{2.00}{33.00}{1}{4.00}{35.00}{2}
\emline{2.00}{33.00}{3}{0.00}{35.00}{4}
\emline{62.00}{33.00}{5}{64.00}{35.00}{6}
\emline{62.00}{33.00}{7}{60.00}{35.00}{8}
\put(2.00,40.00){\makebox(0,0)[cc]{$p$}}
\put(62.00,40.00){\makebox(0,0)[cc]{$G$}}
\put(62.00,0.00){\makebox(0,0)[cc]{$G$}}
\put(2.00,0.00){\makebox(0,0)[cc]{$G$}}
\put(32.00,2.00){\makebox(0,0)[cb]{${\rm Id}$}}
\put(32.00,42.00){\makebox(0,0)[cb]{$e$}}
\put(0.00,20.00){\makebox(0,0)[rc]{$e$}}
\put(64.00,20.00){\makebox(0,0)[lc]{${\rm Id}$}}
\end{picture}
\end{equation}

Now, the diagrams
\begin{equation}
\unitlength=0.50mm
\special{em:linewidth 0.4pt}
\linethickness{0.4pt}
\begin{picture}(221.00,40.00)
\put(8.00,40.00){\vector(1,0){52.00}}
\put(64.00,36.00){\vector(0,-1){32.00}}
\put(4.00,36.00){\vector(0,-1){32.00}}
\put(84.00,40.00){\vector(1,0){52.00}}
\put(84.00,0.00){\vector(1,0){52.00}}
\put(140.00,36.00){\vector(0,-1){32.00}}
\put(80.00,36.00){\vector(0,-1){32.00}}
\put(163.00,40.00){\vector(1,0){52.00}}
\put(163.00,0.00){\vector(1,0){52.00}}
\put(219.00,36.00){\vector(0,-1){32.00}}
\put(159.00,36.00){\vector(0,-1){32.00}}
\put(4.00,40.00){\makebox(0,0)[cc]{$p$}}
\put(0.00,40.00){\makebox(0,0)[rc]{$p\times p=$}}
\put(64.00,40.00){\makebox(0,0)[cc]{$p$}}
\put(64.00,0.00){\makebox(0,0)[cc]{$G$}}
\put(4.00,0.00){\makebox(0,0)[cc]{$G\times G$}}
\put(34.00,2.00){\makebox(0,0)[cb]{$\cdot$}}
\put(62.00,20.00){\makebox(0,0)[rc]{$e$}}
\put(2.00,20.00){\makebox(0,0)[rc]{$e\times e=e$}}
\put(80.00,40.00){\makebox(0,0)[cc]{$p$}}
\put(140.00,40.00){\makebox(0,0)[cc]{$p$}}
\put(140.00,0.00){\makebox(0,0)[cc]{$G$}}
\put(138.00,20.00){\makebox(0,0)[rc]{$e$}}
\put(80.00,0.00){\makebox(0,0)[cc]{$p$}}
\put(110.00,2.00){\makebox(0,0)[cb]{$e$}}
\put(189.00,2.00){\makebox(0,0)[cb]{$i$}}
\put(159.00,0.00){\makebox(0,0)[cc]{$G$}}
\put(219.00,0.00){\makebox(0,0)[cc]{$G$}}
\put(221.00,20.00){\makebox(0,0)[lc]{$e$}}
\put(219.00,40.00){\makebox(0,0)[cc]{$p$}}
\put(159.00,40.00){\makebox(0,0)[cc]{$p$}}
\put(161.00,20.00){\makebox(0,0)[lc]{$e$}}
\put(16.00,0.00){\vector(1,0){44.00}}
\end{picture}
\end{equation}
are commutative, defining thus some arrows
\sar\wt{G}\times\wt{G},\wt{\cdot},\wt G;,
\sar{\bf p},\wt e,\wt G;
and~\sar\wt G,\wt i,\wt G;.

One can easily check that these arrows turn $\wt G$ into a {\sl local}
group, that any morphism $f\colon G\arr G'$ of groups in \cc{C} generates
the morphism $\wt{f}\colon\wt{G}\arr\wt{G'}$ of local
groups and that the correspondence $G\mapsto\wt G$,
$f\mapsto\wt f$ is, in fact, a functor
\begin{equation}
{\rm Loc}\colon \grp C\arr\lgrp{C}\harr
\Groups(\germ C),
\end{equation}
where $\lgrp{C}$ is a full subcategory of
$\Groups(\germ C)$ consisting just of all local groups.
The functor Loc respects finite products.

Besides, for every group $G$ of \cc{C} the arrow (\ref{locincl}) is a morphism
of groups; moreover, the family of all $\iota_G$, where $G$ runs over all
groups of \cc{C}, forms a functor morphism
\begin{equation}
\iota\colon{\rm Loc}\arr U_{{\rm Gr}},
\end{equation}
where the inclusion functor
$U_{{\rm Gr}}\colon\grp C\harr\Groups(\germ C)$
is generated by the inclusion functor $U$ (see (\ref{UandF})) and will be
denoted further simply as $U$, by a standard abuse of notations.

\section{The subgroup $\wt G$ is normal in $G$}
\label{normsub}
Let $G$ and $G'$ be groups in a category \cc{D} and
\sar G',\iota,G; be some morphism of groups,which is mono, considered in
\cc{D}. The arrow $\iota$ will be called a {\bf normal subgroup arrow} if
there exists an arrow $g''\colon G\times G'\arr G'$ in \cc{D} such that
the equality
\begin{equation}
g\cdot\iota g'\cdot g^{-1}=\iota g''
\end{equation}
holds (here $(g,g')={\rm Id}_{G\times G'}$), i.e. the diagram
\begin{equation}
\unitlength=0.50mm
\special{em:linewidth 0.4pt}
\linethickness{0.4pt}
\begin{picture}(96.00,48.00)
\put(4.00,2.00){\makebox(0,0)[cc]{$G'$}}
\put(94.00,2.00){\makebox(0,0)[cc]{$G$}}
\put(94.00,42.00){\makebox(0,0)[cc]{$G\times G\times G$}}
\put(4.00,42.00){\makebox(0,0)[cc]{$G\times G$}}
\put(4.00,38.00){\vector(0,-1){32.00}}
\put(94.00,38.00){\vector(0,-1){32.00}}
\put(16.00,42.00){\vector(1,0){56.00}}
\put(44.00,45.00){\makebox(0,0)[cb]{$(g,\iota g',g^{-1})$}}
\put(51.00,4.00){\makebox(0,0)[cb]{$\iota$}}
\put(2.00,22.00){\makebox(0,0)[rc]{$g''$}}
\put(96.00,22.00){\makebox(0,0)[lc]{$\cdot$}}
\put(8.00,2.00){\vector(1,0){78.00}}
\end{picture}
\end{equation}
is commutative for some $g''$ (unique because $\iota$ is mono). This,
evidently, means that for every object $Z$ of \cc{D} the image
$\iota_*([Z,G']_{\cc{D}})\harr[Z,G]_{\cc{D}}$ is a normal subgroup in
$[Z,G]_{\cc{D}}$.

Note that this is a weak definition of a normal subgroup: the fact that
$\iota$ is a normal subgroup does not imply, generally speaking, the
existence of the ``factorgroup'', because the functor
$Z\mapsto[Z,G]_{\cc{D}}/\iota_*([Z,G']_{\cc{D}})$ need not be
representable.
\begin{Prop}
For every group $G$ in a locus \cc{C} the arrow {\rm (\ref{locincl})} is
a normal subgroup arrow.
\end{Prop}

\section{Actions of groups}
\label{act}
Let \cc{D} be a category with products and \sar M\times G,\cdot,M; be some
arrow of \cc{D}. For every object $Z$ of \cc{D} and a pair of arrows
\sar Z,m,M;, \sar Z,g,G; define the arrow
$m\cdot g\colon Z\arr M$ as the composition arrow
\begin{equation}
\label{mult}
m\cdot g:=(\sar Z,{(m,g)},M\times G;\sar{},\cdot,M;),
\end{equation}
so that, in particular,
\begin{equation}
\cdot=m\cdot g,
\end{equation}
for $(m,g)={\rm Id}_{M\times G}$.

Let now $G$ be a group in \cc{D}. Then the above arrow $\cdot$ is called
a {\bf right action} of the group $G$ on the object $M$ if the following
identities
\begin{equation}
(m\cdot g_1)\cdot g_2=m\cdot(g_1\cdot g_2),\ \ m\cdot1=m
\end{equation}
hold, where $(m,g_1,g_2)={\rm Id}_{M\times G\times G}$ and $m={\rm Id}_M$.

Only right actions of groups will be considered in what follows, so that
the word ``right'' will be omitted.

\section{Local actions of local groups}
\label{locact}
Let \sar M\times G,\cdot,M; be an action of a group $G$ of a locus \cc{C}
on an object $M$ of \cc{C}. Then the normal subgroup arrow (\ref{locincl})
generates the action of the local group $\wt G$ on $M$ as follows:
\begin{equation}
\label{lact}
\sar M\times \wt G,{\rm Id}\times\iota_G,M\times G;\sar{},\cdot,M;.
\end{equation}

Moreover, let ${\bf Z}\colon Z\arr M$ be any object in
\germ C over $M$ (not necessarily open); {\bf Z} can be thought
of as a subobject of the object $M$ of \germ C, the canonical
mono
\begin{equation}
{\bf z}\colon{\bf Z}\arr M
\end{equation}
being defined by the (equivalence class of the) commutative diagram
\begin{equation}
\label{canmon}
\unitlength=0.50mm
\special{em:linewidth 0.4pt}
\linethickness{0.4pt}
\begin{picture}(66.00,48.00)
\put(8.00,42.00){\vector(1,0){52.00}}
\put(8.00,2.00){\vector(1,0){52.00}}
\put(64.00,38.00){\vector(0,-1){32.00}}
\put(4.00,38.00){\vector(0,-1){32.00}}
\put(4.00,2.00){\makebox(0,0)[cc]{$M$}}
\put(64.00,2.00){\makebox(0,0)[cc]{$M$}}
\put(64.00,42.00){\makebox(0,0)[cc]{$M$}}
\put(4.00,42.00){\makebox(0,0)[cc]{$Z$}}
\put(34.00,44.00){\makebox(0,0)[cb]{{\bf Z}}}
\put(34.00,4.00){\makebox(0,0)[cb]{{\rm Id}}}
\put(66.00,22.00){\makebox(0,0)[lc]{Id}}
\put(2.00,22.00){\makebox(0,0)[rc]{{\bf Z}}}
\put(19.00,22.00){\vector(1,0){30.00}}
\put(34.00,24.00){\makebox(0,0)[cb]{{\bf z}}}
\end{picture}
\end{equation}

It turns out that {\bf Z} is an {\bf invariant subobject} with respect to
the action (\ref{lact}), in the sense that there exists (the only, because
{\bf z} is mono) the restriction
\sar{\bf Z}\times\wt G,.,{\bf Z}; of the action (\ref{lact}) such that
the diagram
\begin{equation}
\unitlength=0.50mm
\special{em:linewidth 0.4pt}
\linethickness{0.4pt}
\begin{picture}(72.00,44.00)
\put(70.00,38.00){\vector(0,-1){32.00}}
\put(10.00,38.00){\vector(0,-1){32.00}}
\put(10.00,2.00){\makebox(0,0)[cc]{$M\times\wt G$}}
\put(70.00,2.00){\makebox(0,0)[cc]{$M$}}
\put(70.00,42.00){\makebox(0,0)[cc]{{\bf Z}}}
\put(10.00,42.00){\makebox(0,0)[cc]{${\bf Z}\times\wt G$}}
\put(72.00,22.00){\makebox(0,0)[lc]{{\bf z}}}
\put(8.00,22.00){\makebox(0,0)[rc]{${\bf z}\times{\rm Id}$}}
\put(40.00,44.00){\makebox(0,0)[cb]{.}}
\put(22.00,2.00){\vector(1,0){44.00}}
\put(22.00,42.00){\vector(1,0){44.00}}
\end{picture}
\end{equation}
is commutative.

\section{Difference morphism}
\label{difmor}
Let $G$ be a group in a locus \cc{C} and ${\bf U}\colon U\arr G$ be some
mono. Then, as was explained in the preceding section, the local group
$\wt G$ acts on {\bf U} from the right. Let $G'$ be another group
in \cc{C} and $f\colon{\bf U}\arr G'$ be some arrow (in \germ C).
Define the {\bf difference arrow} (or, rather, left difference arrow) of $f$
\begin{equation}
\Delta f\colon{\bf U}\times\wt G\arr G'
\end{equation}
by means of the expression
\begin{equation}
\Delta f(u,g):=(fu)^{-1}\cdot f(u\cdot g),
\end{equation}
where $(u,g)={\rm Id}_{{\bf U}\times\wt G}$ and the multiplication
$u\cdot g$ etc. is defined by eq.~(\ref{mult}).
\begin{Prop}
\label{extdefdif}
For any pair \sar{\bf Z},u',{\bf U}; and \sar{\bf Z},g',\wt G;
of arrows of \germ C the identity
\begin{equation}
\Delta f(u',g')=(fu')^{-1}\cdot f(u'\cdot g')
\end{equation}
holds.
\end{Prop}
\begin{Cor}
For any $f\colon{\bf U}\arr G'$ and \sar{\bf Z},u',{\bf U}; the identity
\begin{equation}
\Delta f(u',1)=1
\end{equation}
holds. This implies, in turn, that $\Delta f$ uniquely decomposes itself as
follows:
\begin{equation}
\unitlength=0.5mm
\special{em:linewidth 0.4pt}
\linethickness{0.4pt}
\begin{picture}(90.00,46.00)
\put(10.00,2.00){\makebox(0,0)[cc]{${\bf U}\times\wt G$}}
\put(90.00,2.00){\makebox(0,0)[cc]{$G'$}}
\put(50.00,42.00){\makebox(0,0)[cc]{$\wt G'$}}
\put(24.00,2.00){\vector(1,0){61.00}}
\put(14.00,6.00){\vector(1,1){32.00}}
\put(54.00,38.00){\vector(1,-1){32.00}}
\put(50.00,4.00){\makebox(0,0)[lb]{$\Delta f$}}
\put(28.00,24.00){\makebox(0,0)[rb]{$\wt\Delta f$}}
\put(72.00,24.00){\makebox(0,0)[lb]{$\iota_G$}}
\end{picture}
\end{equation}
\end{Cor}

Now we will extend a bit the definition of the difference arrow.
Let ${\bf U}\colon U\arr G$ and ${\bf U'}\colon U'\arr G'$ are objects
of \germ C over groups $G$ and $G'$, respectively. For an arrow
$f\colon {\bf U}\arr {\bf U'}$ define the difference arrow
$\Delta f\colon{\bf U}\times\wt G\arr G'$
as the difference arrow of the composition arrow
\sar{\bf U},f,{\bf U'};\sar{},{\bf u'},G';, where the canonical mono
${\bf u'}$ is defined by the diagram (\ref{canmon}) (with necessary change of
notations there). Similarly, one can define the arrow
$\wt\Delta f\colon U\times\wt G\arr\wt G'$ for this case.

\begin{Prop}
Let $G$, $G'$ and $G''$ be groups of a locus \cc{C}. Let
${\bf U}\colon U\arr G$, ${\bf U'}\colon U'\arr G'$ and
${\bf U''}\colon U''\arr G''$ are objects of \germ C over
$G$, $G'$ and $G''$, respectively. Then for any pair of arrows
$f\colon{\bf U}\arr {\bf U'}$ and $f'\colon{\bf U'}\arr {\bf U''}$
the following identity
\begin{equation}
\wt\Delta(f'f)(u,g)=\wt\Delta f'(fu,\wt\Delta f(u,g))
\end{equation}
holds, where $(u,g)={\rm Id}_{{\bf U}\times G}$.
\end{Prop}
The latter proposition states the chain rule for the difference arrow. It,
clearly, implies that the correspondence
\[({\bf U}\colon U\arr G)\mapsto{\bf U}\times\wt G\]
continues to a functor from the category of ``subobjects over groups'' of
\cc{C} to the category of pointed bundles (see sect. \ref{gr}) of the category
\germ C.

For any integer $n>1$ the $n$-{\bf th difference arrows}
$\Delta^nf$ and $\wt\Delta^nf$ can be defined recurrently as follows:
\begin{userlabel}
\begin{eqnarray}
\label{n-thdiff}
\Delta^nf=\Delta(\Delta^{n-1}f)\colon\U\times\wt G^n\arr G'\,,\\
\wt\Delta^nf=\wt\Delta(\wt\Delta^{n-1}f)\colon\U\times\wt G^n\arr \wt G'\,.
\end{eqnarray}
\end{userlabel}

\section{Partial difference morphisms}
\label{partdiff}
Consider now the situation when ${\bf U}$ is of the form
${\bf U}\colon U\arr G_1\times\dots\times G_n$, where $G_1$,\dots,$G_n$ are
some groups in \cc{C}. There are canonical morphisms of groups
\begin{equation}
\label{caningr}
i_\alpha\colon\wt G_\alpha\arr\wt G_1\times\dots\times
\wt G_n\ \ (\alpha=1,\dots,n)
\end{equation}
defined as $i_\alpha(x)=(1,\dots,x,\dots,1)$,
where $x={\rm Id}_{G_\alpha}$.

Now, given ${\bf U'}\colon U'\arr G'$ and an arrow
$f\colon{\bf U}\arr{\bf U'}$ define the $\alpha$'s {\bf partial difference
morphism}
\begin{equation}
\Delta_\alpha f\colon{\bf U}\times\wt G_\alpha\arr G'
\end{equation}
as follows:
\begin{equation}
\Delta_\alpha f=\Delta f(u,i_\alpha g_\alpha)\ \
(\alpha=1,\dots,n),
\end{equation}
where $(u,g_\alpha)={\rm Id}_{{\bf U}\times\wt G_\alpha}$.

For arbitrary arrows \sar{\bf Z},u',{\bf U}; and \sar{\bf Z},g',\wt G;
one has the counterpart of the equality (\ref{extdefdif}):
\begin{equation}
\Delta_\alpha f(u',g')=
(fu')^{-1}\cdot f(u'\cdot i_\alpha g')\ \ (\alpha=1,\dots,n).
\end{equation}

Besides, one can check that $\Delta f$ can be expressed in terms of
partial difference arrows $\Delta_\alpha f$ as follows:
\begin{equation}
\label{difptdif}
\begin{tabular}{lll}
$\Delta f$&=&$\Delta_1f(u,g_1)\cdot\dots
\cdot\Delta_\alpha f(u\cdot i_1g_1\cdot\dots
\cdot i_{\alpha-1}g_{\alpha-1},g_\alpha)\cdot\dots$\\
&&$\dots\cdot\Delta_nf(u\cdot i_1 g_1\cdot\dots\cdot i_{n-1}g_{n-1},g_n)$,
\end{tabular}
\end{equation}
where
$(u,g_1,\dots,g_n)={\rm Id}_{{\bf U}\times\wt G_1\times\dots\times
\wt G_n}$ and arrows $i_\alpha$ are defined above
(see (\ref{caningr})).

\section{Local differential calculus}
\label{calc}
We can extract, at last, the archetype of ``local approximation of maps
by linear maps'' formalizing it in the general notion of {\sl differential
calculus}. First, it is clear now that the thing that is to be ``locally''
approximated is not an arrow $f\colon{\bf U}\arr{\bf U'}$ itself but, rather,
its difference arrow $\wt\Delta f$. On the other hand, we are,
generally speaking, to deal with groups equipped with some additional
structure, instead of ``pure'' groups of the locus \cc{C} (in an archetypical
example of ``the'' calculus this additional structure is the structure of
a finite-dimensional vector space).

Now the definitions.
Let \cc{V} be a category with finite products and zero object
(i.e. both terminal and initial).
Call a {\bf precalculus} on the locus \cc{C}
any structure functor
$G\colon\cc{V}\arr\grp C$
respecting finite products.
Recall that a functor is called a structure functor if it is
fully faithful and, besides, there exist direct and inverse images of
structures along iso's (see, e.g., \cite{Man}). Given a precalculus $G$
on \cc{C}
define the category ${\rm C}^0(\cc{C},G)$ as follows. Objects of
${\rm C}^0(\cc{C},G)$ are all pairs $\la V,{\bf U}\ra$ such that
${\bf U}\colon U\arr GV$ is an object of \germ C over the group
$GV$ or, rather, over its underlying object in \cc{C}; we will often omit in
what follows the forgetful functor $G$. An arrow
$f\colon\la V,{\bf U}\ra\arr\la V',{\bf U'}\ra$ in
${\rm C}^0(\cc{C},G)$ is simply any
triple $\la V,\ul f,V'\ra$, where $\ul f$ is an
arrow (in \germ C)
between {\bf U} and ${\bf U'}$. In accord with general math practice,
we would often write $f$ instead of $\ul f$ and vice versa,
hoping that it would not lead to confusion.

The full subcategory of the category ${\rm C}^0(\cc{C},G)$ consisting
of all objects ${\bf U}\colon U\arr V$ which are open arrows of \cc{C} will
be called the {\bf category of open regions of class} ${\rm C}^0$
of the precalculus $G\colon\cc{V}\arr \cc{C}$ and
will be denoted as ${\rm Reg}^0(\cc{C},G)$.

The set-valued functor
\begin{equation}
\cc{H}\colon{\rm C}^0(\cc{C},G)^{op}\times\cc{V}^{op}\times\cc{V}\arr{\bf Set}
\end{equation}
sending a triple $\la{\bf U},V,V'\ra$ into the hom-set
$[{\bf U}\times\wt{GV},\wt{GV'}]$
has the evident natural structure of the group in the corresponding
functor category. For any arrow $f\colon{\bf U}\arr{\bf U'}$ in
${\rm C}^0(\cc{C},G)$ between objects ${\bf U}\colon U\arr V$ and
${\bf U'}\colon U'\arr V'$ the (local) difference arrow $\wt\Delta f$
clearly belongs to $\cc{H}({\bf U},V,V')$.

A precalculus $G\colon\cc{V}\arr \grp C$ together with a
pair $\la\wt{\cc{L}},o\ra$ of subfunctors of \cc{H} will be
called a {\bf local (differential) calculus (on \cc{C})}
if $o$ is a subgroup of \cc{H} such that
the following conditions are satisfied:

\noindent{\bf(DC1)} The subfunctor $\wt{\cc{L}}$ contains the unity
1 of the group \cc{H}.

\noindent{\bf(DC2)} The intersection $\la\wt{\cc{L}}\ra\cap o$
is
a trivial subgroup of \cc{H}. Here $\la\wt{\cc{L}}\ra$ denotes
the subgroup of \cc{H} {\sl generated} by the subfunctor $\wt{\cc{L}}$;

\noindent{\bf(DC3)} If
$l\colon{\bf U}\times \wt V\arr \wt{V'}$
belongs to
$\wt{\cc{L}}({\bf U},V,V')$
and
$l'\colon{\bf U}\times \wt{V'}\arr \wt{V''}$
belongs
to $\wt{\cc{L}}({\bf U},V',V'')$,
then the composition arrow
$l'(u,l)$
belongs to
$\wt{\cc{L}}({\bf U},V,V'')$
(here $u={\rm Id}_{\bf U}$);

\noindent{\bf(DC4)} If $l$, $l'$ and $u$ are as above and
$i\colon{\bf U}\times\wt V\arr\wt{V'}$ belongs to $o({\bf U},V,V')$ then
\begin{equation}
l'(u,l\cdot i)=l'(u,l)\cdot i'
\end{equation}
for some $i'\in o({\bf U},V,V'')$;

\noindent{\bf(DC5)} If $l$ and $i$ are as above and
$i'\colon{\bf U}\times\wt{V'}\arr\wt{V''}$
belongs to
$o({\bf U},V',V'')$ then the
composition arrow $i'(u,l\cdot i)$ belongs to
$o({\bf U},V,V'')$;

\noindent{\bf(DC6)} The subgroup $o$ is normal with respect to
the subfunctor \cc{L} in the sense that for any
$l$ and $i$ as above there exists a $j\in o({\bf U},V,V')$ such that
the equality $i\cdot l=l\cdot j$ holds;

\noindent{\bf(DC7)} For any object {\bf U} of ${\rm C}^0(\cc{C},\cc{V})$
and any object $V$ of \cc{V} the projection arrow
\sar{\bf U}\times\wt V,g,\wt V; belongs to $\wt{\cc{L}}({\bf U},V,V)$;

\noindent{\bf(DC8)} Both $\wt{\cc{L}}$ and $o$ are closed with respect to
products: if $l_a\colon{\bf U}\times\wt V\arr\wt{V_a}$ belongs to
$\wt{\cc{L}}({\bf U},V,V_a)$ for $a=1,2$ then the arrow
$(l_1,l_2)\colon {\bf U}\times\wt V\arr\wt{V_1}\times\wt{V_2}$ belongs to
$\wt{\cc{L}}({\bf U},V,V_1\times V_2)$.

Note, that we do not require in the definition above neither that
$\wt{\cc{L}}$ should be a subgroup of \cc{H}, nor that families belonging to
\cc{L} should be families of group morphisms: both this requirements would be
incompatible, generally speaking, with each other, excepting the case when
the precalculus $G\colon\cc{V}\arr \grp C$
pulls through the category \agrp{C} of abelian groups of \cc{C};
we will call such a precalculus $G\colon\cc{V}\arr\grp C$ {\bf additive} if,
besides, the category \cc{V} is additive as well as the restriction of the
functor $G$ on the subcategory \agrp{C} of \grp C.
Call a calculus $\la G,\wt{\cc{L}},o\ra$
{\bf multiplicative}, resp. {\bf quasimultiplicative} if for any triple
${\bf U}$,\cc{V},\cc{V'} and any
$l\in\wt{\cc{L}}({\bf U},\cc{V},\cc{V'})$ the family $l$
is a family of group morphisms, resp. has a representation
\begin{equation}
l=l_1\cdot\dots\cdot l_n\,,
\end{equation}
such that every $l_\alpha$ is either a family of group morphisms or inverse
to such a family.

Of course, quasimultiplicative calculus is multiplicative for the case of
additive calculi.

It is clear as well that $\wt{\cc{L}}$, resp. $o$ define due to (DC3),
resp. due to (DC5) a subcategory in the category of pointed trivial bundles
in the
category ${\rm C}^0(\cc{C},G)$: namely, of all arrows between trivial bundles
with principal part belonging to $\wt{\cc{L}}$, resp. to $o$. Such
bundles will be called $\wt{\cc{L}}${\bf-bundles},
resp. $o${\bf-bundles}.

An arrow $f\colon{\bf U}\arr{\bf U'}$ of ${\rm C}^0(\cc{C},G)$ will be called
{\bf locally differentiable} or {\bf locally of class}
${\rm C}^1$ (in the calculus
$\la G,\wt{\cc{L}},o\ra$) if its local difference arrow
$\wt{\Delta}f$ presents itself as
\begin{equation}
\label{difder}
\wt{\Delta}f=\wt Df\cdot i\ \
(\wt Df\in\wt{\cc{L}},\ i\in o).
\end{equation}
It then follows from (DC2) that this decomposition is unique; the family
$\wt Df$ will be called the {\bf local derivative} of $f$.

\begin{Prop}
\label{diffcomp}
Let ${\bf U}\colon U\arr V$ be an object of ${\rm C}^0(\cc{C},G)$  over $V$.
If $f\colon{\bf U}\arr{\bf U'}$ and $f'\colon{\bf U'}\arr{\bf U''}$ are
locally differentiable then their composition $f'f$ is locally differentiable
and the chain rule
\begin{userlabel}
\begin{equation}
\wt D(f'f)=\wt Df'(fu,\wt Df)\ \
((u,v)={\rm Id}_{{\bf U}\times V})
\end{equation}
is valid.

More generally, for any arrows \sar Z,u,\U; and \sar Z,v,V; one has
the identity
\begin{equation}
\label{}
\wt D(f'f)(u,v)=\wt Df'(fu,\wt Df(u,v))\,.
\end{equation}
\end{userlabel}
\end{Prop}

Prop.\ref{diffcomp} implies that locally differentiable arrows
form a subcategory of ${\rm C}^0(\cc{C},G)$ with the same (due to (DC7))
set of objects as ${\rm C}^0(\cc{C},G)$ itself; this subcategory
will be denoted ${\rm C}^1(\cc{C},G,\wt{\cc{L}},o)$. Moreover,
the correspondence
\begin{equation}
(f\colon{\bf U}\arr{\bf U'})\mapsto(fu,\wt Df),
\end{equation}
where {\bf U} is an object over $V$ and $u\colon{\bf U}\times V\arr{\bf U}$ is
the projection, continues to the functor $\wt T$ (the {\bf local tangent functor})
from the category
${\rm C}^1(\cc{C},G,\wt{\cc{L}},o)$ to the category of trivial
$\wt{\cc{L}}$-bundles of the category ${\rm C}^0(\cc{C},G)$.

Let {\bf U} be an object of ${\rm C}^0(\cc{C},G)$ over a finite product
$V_1\times\dots\times V_n$ of objects of \cc{V} and
$f\colon{\bf U}\arr{\bf U'}$ be an arrow of ${\rm C}^0(\cc{C},G)$. Call $f$
{\bf locally differentiable along} $V_\alpha$ ($\alpha=1,\dots,n$),
if the local partial difference $\wt{\Delta}_\alpha f$ presents itself as
\begin{equation}
\wt{\Delta}_\alpha f=\wt D_\alpha f\cdot i\ \
(\wt D_\alpha f\in\wt{\cc{L}},\ i\in o).
\end{equation}
The family $\wt D_\alpha f$ will be called the {\bf local partial
derivative} of $f$ along $V_\alpha$.

The following proposition states that for {\sl multiplicative} calculi
local differentiability implies the existence of local partial derivatives.
\begin{Prop}
Let $\la G,\wt{\cc{L}},o\ra$ be a multiplicative calculus.
Then:

\noindent{\rm(a)} If {\bf U} is an object of ${\bf C}^0(\cc{C},G)$ over the
finite product
$V_1\times\dots\times V_n$ of objects of \cc{V} and
$f\colon{\bf U}\arr{\bf U'}$ is locally differentiable then for any
$\alpha=1,\dots,n$
there exists the local partial derivative $\wt D_\alpha f$ of $f$ along
$V_\alpha$ and the identity
\begin{equation}
\wt D_\alpha f=\wt Df(u,i_\alpha v_\alpha)
\end{equation}
holds, where
$(u,v_\alpha)={\rm Id}_{{\bf U}\times\wt V_\alpha}$
and
$i_\alpha\colon\wt V_\alpha\arr\wt V_1\times\dots\times
\wt V_n\ \ (\alpha=1,\dots,n)$
is the canonical injection arrow 
$i_\alpha(x)=(1,\dots,x,\dots,1)$ for $x={\rm Id}_{V_\alpha}$%
.

\noindent{\rm(b)} Besides, the local derivative $\wt Df$ is expressed
through local partial derivatives as follows:
\begin{equation}
\wt Df=\wt D_1f(u,v_1)\cdot\dots\cdot\wt D_nf(u,v_n),
\end{equation}
where $(u,v_1,\dots,v_n)={\rm Id}_{{\bf U}\times V_1\times\dots\times V_n}$.

\noindent{\rm c)}
For any $\alpha,\beta=1,\dots,n$ the local partial derivatives
$\wt D_\alpha f$ and $\wt D_\beta f$ commute:
\begin{equation}
\wt D_\alpha f(u,v_\alpha)\cdot\wt D_\beta f(u,v_\beta)=
\wt D_\beta f(u,v_\beta)\cdot\wt D_\alpha f(u,v_\alpha).
\end{equation}
\end{Prop}

\section{Arrows of class $C^n$}

What about arrows of class $C^n$ with $n=2,\dots,\infty$?

One of the possible definitions is the classical recurrent definition:
for an integer $n>1$ we say that an arrow
$f\colon{\bf U}\arr{\bf U'}$ of ${\rm C}^0(\cc{C},G)$
is {\bf locally of class}
${\rm C}^n$ or is $n$ {\bf times locally differentiable}
(in the calculus
$\la G,\wt{\cc{L}},o\ra$), if it is localy of class $C^1$
and its local derivative arrow $\wt Df$ is of class $C^{n-1}$;
$f$ is {\bf of class $C^\infty$} if it is of class $C^n$ for any $n\ge1$.

This defines, for any $n=1,2,\dots,\infty$, the category
$C^n(\cc{C},G,\wt{\cc{L}},o)$
of arrows $n$ times locally differentiable
in the calculus $\la G,\wt{\cc{L}},o\ra$.

\section{Global differential calculus}

Above definitions of local differential calculus are too broad:
the categories
$C^n(\cc{C},G,\wt{\cc{L}},o)$
contain many ``non-classical'' objects, in fact, all germs of monos
over the category $\cc{V}$.

For example, in an archetypical case
of $\cc{V}$ consisting of all finite-dimensional vector spaces
over real or comlex numbers, these include such ``model objects''
as regions with boundaries or angles, as well as closed subsets
of an open region, giving rise to the corresponding theory of
smoothness.

If we want to restrict ourselves to the ``classical'' subcategory
${\rm Reg}^n(\cc{C},G,\wt{\cc{L}},o)$
of {\bf open regions with arrows locally of class} $C^n$,
we would like to globalize as well local derivatives
or, which is the same, the local tangent functor $\wt T$.

This means that for any locally differentiable arrow
\sar\U,f,V'; with $\U=\sar U,u,GV;$ we would like to have
the uniquelly determined lifting
$Df\colon\U\times V\arr V'$
of the local derivative arrow
$\wt Df\colon\U\times \wt V\arr \wt{V'}$.

An evident way to do this is as follows.

For a precalculus
$G\colon\cc{V}\arr\grp C$
let $H$ be the functor
\begin{equation}
H\colon{\rm C}^0(\cc{C},G)^{op}\times\cc{V}^{op}\times\cc{V}\arr{\bf Set}
\end{equation}
sending a triple
$\la\U,V,V'\ra$
into the subset of the hom-set
$[\U\times GV,V']$
consisting of all arrows
\sar\U,f,V; such that $f(u,1)=1$ for any \sar Z,u,\U;.
The functor $H$ has an evident structure of the group
in the corresponding functor category.

Call a subfunctor \cc{L} of the functor $H$ {\bf suitable}
if it satisfies the following counterpart of conditions
(DC1), (DC3), (DC7) and (DC8) above:

\noindent{\bf(GDC1)} \cc{L} contains the unity $1$ of the group $H$;

\noindent{\bf(GDC3)} \cc{C} is stable under composition:
if
$l\colon{\bf U}\times  V\arr {V'}$
belongs to
$\cc L({\bf U},V,V')$
and
$l'\colon{\bf U}\times {V'}\arr {V''}$
belongs
to $\cc L({\bf U},V',V'')$,
then the composition arrow
$l'(u,l)$
belongs to
$\cc L({\bf U},V,V'')$
(here $u={\rm Id}_{\bf U}$).

\noindent{\bf(GDC7)} For any object {\bf U} of ${\rm C}^0(\cc{C},\cc{V})$
and any object $V$ of \cc{V} the projection arrow
\sar{\bf U}\times V,g, V; belongs to $\cc{L}({\bf U},V,V)$;

\noindent{\bf(GDC8)} \cc{L} is closed with respect to
products: if $l_a\colon{\bf U}\times V\arr V_a$ belongs to
$\cc{L}({\bf U},V,V_a)$ for $a=1,2$ then the arrow
$(l_1,l_2)\colon {\bf U}\times V\arr V_1\times V_2$ belongs to
$\cc{L}({\bf U},V,V_1\times V_2)$. The similar is true for $o$.

It is clear that any suitable subfunctor \cc L of $H$ defines
a subcategory in the category of pointed trivial bundles
in the
category ${\rm C}^0(\cc{C},G)$: namely, of all arrows between trivial bundles
with principal part belonging to \cc L. Such
bundles will be called \cc L{\bf-bundles}.

A suitable subfunctor \cc L of $H$ canonically generates, besides,
a subfunctor \wt{\cc L} by restricting any arrow
\sar\U\times V,l,V'; to \sar\U\times\wt V,1\times\iota,\U\times\wt V;;
this restriction pulls then through the unique arrow
\sar\U\times\wt V,\wt l,\wt{V'}; belonging to \cc{H},
as follows from the definition of the functor $H$.

A precalculus
$G\colon\cc{V}\arr\grp C$
together with a pair of subfunctors
$\la L\subset H,o\subset\cc{H}\ra$
will be called
a {\bf(global) calculus}
if the subfunctor $L\subset H$ is suitable
and the pair
$\la \wt L\subset\cc{H},o\subset\cc{H}\ra$
is a local calculus.

\section{The first example --- at last!}

The first example of a (global) differential multiplicative
calculus, being rather trivial, exists in every category \cc{C}.
In fact, it is functorial on \cc{C}.

An arrow \sar G,f,G'; between
groups of the category \cc{C} will be called
a {\bf left affine arrow}
iff there exist an element \sar p,c,G'; of the group $G'$ (with $p$
being the final object of \cc{C}) and a morphism of groups
\sar G,f,G'; such that $f=c\cdot f'$.
{\bf Right affine arrows} are defined similarly.

One easily sees that $f$ is a left affine arrow iff
it is a right affine arrow, so that we will call such
arrows simply {\bf affine}.

\begin{Prop}
Let \cc{C} be any category.
Let \cc{O} consists of all isomorphisms of
the category \cc{C}.
Let \cc{V} be the category of all groups of the category \cc{C},
and \sar\cc{V},G,\grp C; be an identity functor.
Let the subgroup $o$ of $H$
be the trivial subgroup of $H$.
Let the subfunctor \cc{L} of $H$ consists of all
families of group morphisms. Then:

\noindent{\rm a)} The functor $G$ together with the pair
$\la\cc{L},o\ra$ is a global differential calculus.

\noindent{\rm b)} The category
${\rm Reg}^1(\cc{C},G,\cc{L},o)$ of regions of class $C^1$
has the following description:
an arrow \sar G,f,G'; between
groups of the category \cc{C} is of class $C^1$,
iff it is affine.

\noindent{\rm c)} Besides, any arrow of class $C^1$ is of
class $C^\infty$.
\end{Prop}


\section{Example 2: classical differential calculus}

We suggest that the reader who was patient enough to read the text up to this point
can reinterprete without difficulties the ``ordinary'' calculus (finite-dimensional,
in Banach spaces, in locally convex spaces, etc.) in terms of an appropriate functor $G$
and construct an adequate differential calculus $\la\cc{L},o\ra$ for all of these cases :).

\end{document}